\documentclass[12 pt]{amsart}

\usepackage{amscd,amsfonts,amssymb,amsmath}
\usepackage{tikz-cd}
\usepackage{xcolor}
\usepackage{placeins}
\usepackage{hyperref}

\newtheorem{theorem}{Theorem}[section]
\theoremstyle{definition}
\newtheorem{conjecture}[theorem]{Conjecture}

\newtheorem{proposition}[theorem]{Proposition}
\newtheorem{question}[theorem]{Question}

\newtheorem{definition}[theorem]{Definition}
\newtheorem{example}[theorem]{Example}

\numberwithin{equation}{section}

\newcommand{\im}{\operatorname{Im}}

\newcommand{\id}{\mathrm{id}}

\title{Algebra, group, and Hopf Rota--Baxter operators
}
\author[Valeriy G.~Bardakov,  Igor M.~Nikonov, and Viktor N. Zhelaybin]{Valeriy G.~Bardakov,  Igor M.~Nikonov, and Viktor N. Zhelaybin}

\date{\today}

\begin{document}
\sloppy



\maketitle
\begin{abstract} We  know definition of  Rota--Baxter operators on  different algebraic systems. For examples, on groups, on algebras,  on Hopf algebras.
On some  algebraic systems it is possible to define different types of Rota--Baxter operators. For example, on group algebra it is possible to define Rota--Baxter operator as on associative algebra, group Rota--Baxter operator and Rota--Baxter operator as on a Hopf algebra. We are investigating the following question: What  are connections between these operators?
We are studying these questions for the Sweedler algebra $H_4$, that is 4-dimension non--cocommutative Hopf algebra.


\end{abstract}

\maketitle
\tableofcontents

\section{Introduction}

Rota--Baxter operators on algebras were introduced in  the middle of the previous century~\cite{Baxter,Tricomi}. These operators  have different  connections
with mathematical physics, Hopf algebras,
 operad theory, and with some other branches of mathematics,  see the monograph~\cite{GuoMonograph}.
Rota--Baxter operators on groups were  introduced in 2020 by L.~Guo, H.~Lang, Y.~Sheng~\cite{Guo2020}.
After this  work, the study of Rota--Baxter operators on groups  continued in~\cite{BG2, BG1, Goncharov2020,JSZ}.
In particular, M.~Goncharov \cite{Goncharov2020} introduced and studied a Rota--Baxter operator on a cocommutative Hopf algebra
that generalizes the notions of a Rota--Baxter operator on a group and a Rota--Baxter operator of
weight 1 on a Lie algebra.
Further relative Rota--Baxter operators on groups were defined \cite{JSZ}, which are generalization of Rota--Baxter operators on groups and
 their connections with relative Rota--Baxter operators on algebras were found. In \cite{BG2} a connection between Rota--Baxter groups and skew left braces was found.
Similar result for relative Rota--Baxter groups was found in \cite{BGST, RS-23}.
In \cite{BN} relative Rota--Baxter operators and Rota--Baxter operators on arbitrary Hopf algebra were defined.
Averaging and relative Rota--Baxter  operators on racks, quandles  and rack algebras were defined in \cite{BB}.

We know  definitions of  RB--operators on different algebraic systems: on algebras, on groups, on group algebras, on Hopf algebras.
It is known \cite{Guo2020} that  a group RB--operator on Lie group induces an algebra RB--operator on the corresponding Lie algebra.
Another example can be found in  \cite{BG2}, where for a RB-group $(G;B)$,  a Lie algebra $L(G)$ on the quotients of the lower central series was constructed
and the authors proved that $B$ induces an algebra RB--operator on $L(G)$.
On some algebras it is possible to define different types of Rota--Baxter operators. For example, if $G$ is a group, $\Bbbk[G]$ is its group algebra, then we can extend by linearity a group RB--operator from $G$ to $\Bbbk[G]$. On the other hand side, we can consider   $\Bbbk[G]$ as a cocommutative Hopf algebra and find a  RB--operator  using construction from \cite{Goncharov2020}.

The next questions also seem interesting. Let $A$ be an associative algebra. We can find its associative algebra RB--operators (AA--RB--operators). We will call these operators  {\it associative algebra RB--operators}. Further, we can construct the Lie algebra $A^{(-)}$ and find its RB--operators, which we will call  {\it Lie  algebra RB--operators} (LA--RB--operators). Also, we can construct the Jordan algebra  $A^{(+)}$ and find its RB--operators, which we will call  {\it Jordan algebra RB--operators} (JA--RB--operators).  It is interesting to find connections between them. It is easy to prove that any associative algebra RB--operator on $A$ is a Lie algebra RB--operator on $A^{(-)}$ and is a Jordan algebra RB-operator on $A^{(+)}$.

Let us give a short outline of the paper.

In~Section~\ref{Prel}, we state the required preliminaries on Hopf algebras, Rota--Baxter operators on algebras, on  groups, and on Hopf algebras. We prove
in Proposition \ref{p2.2} that if $q^2 = -1$, then there is a non-trivial homomorphism $U_q (sl_2) \to H_4$ of the Hopf algebras.

In Section \ref{G-A-RBO}, we   prove (see Theorem \ref{GA-oper}) that if $B \colon G \to G$ is a group RB--operator on a group $G$ and its linear extension to  a group algebra $\Bbbk[G]$ is an algebra RB--operator, then
 $G=KH$ is an exact factorization;
$B(kh)=h$ for all $k \in K$, $h\in H$, $H$
 is commutative and $h^2Kh^{-2}\subset K$ for all $h\in H$.
Conversely, for any exact group factorization $G=KH$ such that $H$ is commutative and $h^2Kh^{-2}\subset K$ for all $h\in H$, the operator $B(kh)=h$ is a group
RB--operator on $G$ and its linear extension defines an algebra RB--operator on~$\Bbbk[G]$.

In Section \ref{RB_H4}, we  study  Rota--Baxter operators on Hopf algebra $H_4$ and the corresponding Lie algebra $H_4^{(-)}$. The algebra $H_4$ contains a cyclic group of order 2. On such group there exist two RB--operators $B_0$ and $B_{-1}$. The fist one sends all elements to the unit and the second one is the identity map. We find all extensions of these operators to $H_4$ which give  RB--operators on cocommutative Hopf algebra $H_4$. As a consequence of these results, we prove, in Proposition \ref{notHom}, that there exist an RB--operator on $H_4$ which is not a coalgebra  homomorphism. At the end of this section we check which  RB--operators on $H_4$ are Lie algebra RB--operators on the Lie algebra $H_4^{(-)}$.
The full description of the Lie algebra RB--operators on  $H_4^{(-)}$ can be found in \cite{BB1}.

In last section we formulate some questions for future research.

\bigskip


\section{Preliminaries} \label{Prel}

\subsection{Hopf algebras} Recall some facts about  Hopf algebras (see, for example, \cite[Chapter~11]{T}).

 Let $A$ be an associative algebra with unit $1_A$ over a commutative ring with unit, $\Bbbk$. Assume that $A$ is provided with multiplicative
$\Bbbk$-linear homomorphisms
$$
\Delta \colon A \to A^{\otimes 2} = A \otimes_{\Bbbk} A~\mbox{and}~\varepsilon \colon A \to \Bbbk,
$$
called
the comultiplication and the counit respectively, and a $\Bbbk$-linear map $S \colon A \to A$, called the antipode. It is understood that $\Delta(1_A) = 1_A \otimes 1_A$ and $\varepsilon(1_A) = 1$. The tuple $(A, \Delta, \varepsilon, S)$ is said to be a Hopf  algebra if these
homomorphisms satisfy together with the multiplication $m \colon A \times A \to A$ the following identities:
\begin{equation}\label{1.1.a}
(\id_A  \otimes \Delta) \Delta  = (\Delta \otimes \id_A) \Delta,
\end{equation}
\begin{equation}\label{1.1.b}
m (S \otimes \id_A ) \Delta  = m (\id_A \otimes S) \Delta = \varepsilon \cdot 1_A,
\end{equation}
\begin{equation}\label{1.1.c}
(\varepsilon \otimes \id_A ) \Delta  = (\id_A \otimes  \varepsilon) \Delta = \id_A.
\end{equation}
The first identity means that the comultiplication is coassociative.
Note that to  write down the first equality we identify $(A \otimes A) \otimes A = A \otimes (A \otimes A)$ via $(a \otimes b) \otimes c = a \otimes (b \otimes c)$, where $a, b, c \in A$.  Similarly, to write down the third equality, we identify $\Bbbk \otimes A = A \otimes \Bbbk = A$ via $1 \otimes a = a \otimes 1 = a$. The axioms imply that the antipode $S$ is an antiautomorphism of both the
algebra and the coalgebra structures in~$A$. This means that
$$
m (S \otimes S) = S \circ m \circ P_A \colon A^{\otimes 2} \to A,~~~~P_A (S \otimes S) \Delta = \Delta  \circ S  \colon A \to A^{\otimes 2},
$$
where $P_A$ denotes the flip
$$
a \otimes b \mapsto b \otimes a \colon A^{\otimes 2} \to A^{\otimes 2}.
$$
 It also follows from the axioms that $S(1_A) = 1_A$ and $\varepsilon \circ S = \varepsilon \colon A \to \Bbbk$.

\begin{example}
\begin{enumerate}
\item A group ring $\Bbbk[G]$ is a Hopf algebra if the homomorphisms $\Delta, \varepsilon, S$  are defined on the additive generators $g \in  G$ by the formulas
$$
\Delta(g) = g \otimes g,~~ \varepsilon(g) = 1,~~ S(g) = g^{-1}.
$$

\item The universal enveloping algebra $U(\mathfrak{g})$ of a Lie algebra $\mathfrak{g}$ is a Hopf algebra, if the homomorphisms $\Delta, \varepsilon, S$  are defined on the multiplicative  generators $g \in \mathfrak{g}$ by the formulas
$$
\Delta(g) = g \otimes 1 + 1 \otimes g,~~ \varepsilon(g) = 0,~~ S(g) = -g.
$$
\item The famous non-cocommutative 4-dimension Sweedler  algebra $H_4$ as algebra is generated by two elements $x$, $g$ with multiplication
$$
g^2 = 1,~~x^2 = 0,~~x g = -g x.
$$
Comultiplication, counit and antipode are defined by the rules
$$
\Delta(g) = g \otimes g,~~\Delta(x) = x \otimes 1 + g \otimes x,~~ \Delta(g x) = g x \otimes g + 1 \otimes g x,
$$
$$
\varepsilon(g) = 1,~~ \varepsilon(x) = 0,~~ \varepsilon(gx) = 0,
$$
$$
 S(g) = g^{-1} = g,~~S(x) = -g x ,~~ S(g x) = x .
$$

The antipode $S$ has order 4 and for any $a \in H_4$ we have $S^2(a) = g a g^{-1}$.
\end{enumerate}
\end{example}

 Throughout the paper  we shall use the
Sweedler's notation with superscripts for all the comultiplications,
$$
\Delta(h) = h_{(1)} \otimes h_{(2)},~~~h \in H.
$$

\medskip

Let us fix an invertible element $q$ of $\Bbbk$ different from $1$ and $-1$ so that the fraction $\frac{1}{q - q^{-1}}$ is well-defined.
Recall (see, for example \cite[Capters 6-7]{Kas}) that the Hopf algebra  $U_q(sl_2)$ as algebra is generated by four elements $E$, $F$, $K$, $K^{-1}$ with the relations:
$$
K K^{-1} = K^{-1} K = 1,
$$
$$
K E K^{-1} = q^2 E,~~K F K^{-1} = q^{-2} F,
$$
$$
[E, F] = \frac{K - K^{-1}}{q - q^{-1}},
$$
$$\Delta(E)=1\otimes E+E\otimes K, ~~~\Delta(F)=K^{-1}\otimes F+F\otimes 1,$$
$$\Delta(K)=K\otimes K, ~~~\Delta(K^{-1})=K^{-1}\otimes K^{-1},$$
$$\epsilon(E)=\epsilon(F)=0, ~~~\epsilon(K)=\epsilon(K^{-1})=1,$$
$$S(E)=-EK^{-1}, ~~S(F)=-KF, ~~S(K)=K^{-1}, ~~S(K^{-1})=K.$$

Let $I$ be an ideal of  $U_q(sl_2)$ generated by $K^2-1$. Then $I$ is a coideal of Hopf algebra $U_q(sl_2).$ Indeed,
$$\Delta(K^2-1)=K^2\otimes K^2-1\otimes 1=(K^2-1)\otimes 1+1\otimes (K^2-1),$$
$$S(K^2-1)=-K^{-2}(K^2-1).$$ Therefore, $I$ is a coideal of Hopf algebra $U_q(sl_2).$

\begin{proposition} \label{p2.2}
If $q^2 = -1$, then the map $\varphi \colon U_q (sl_2) \to H_4$ that is defined on the generators
$$
1 \mapsto 1,~~K \mapsto g,~~E \mapsto gx,~~F \mapsto x,
$$
defines a homomorphism of Hopf algebras with $\ker \varphi=I$.
\end{proposition}

\begin{proof}
It is need to check that the  relations of quotient algebra $U_q (sl_2)/I$ go to relations of $H_4$ under $\varphi$.

The relation $[E, F] = 0\, mod\, I $ goes to relation $[gx, x] =0$, Since
$$
[gx, x] =gx x - x gx =0,
$$
this relation hold.

The relation $K E K^{-1} = q^2 E$ goes to relation $g gx g = q^2 g x$, Since $x g = -x g$ and $q^2 = -1$,
this relation hold.

The relation $K F K^{-1} = q^2 F$ goes to relation $g x g = q^2  x$, Since $x g = -x g$ and $q^2 = -1$,
this relation hold.

The relation $\Delta (E) = 1 \otimes E + E  \otimes K$ goes to relation $\Delta (gx) = 1 \otimes gx + gx  \otimes g$, which is true by the definition multiplication on $H_4$.

The relation $\Delta (F) = K^{-1} \otimes F + F  \otimes 1$ goes to relation $\Delta (x) = g \otimes x + x  \otimes 1$, which is true by the definition multiplication on $H_4$.

The counit relations $\varepsilon (E) = \varepsilon (F) = 0$ go to relations $\varepsilon (g x) = \varepsilon (x) = 0$ which are true.

An antipode relation $S(E) = - E K^{-1}$ goes to relation  $S(gx) = - gxg$. Since $- gxg = x$, this relation holds.

The checking for other relations is similar.
\end{proof}

\subsection{Rota--Baxter operators}
Rota--Baxter operators for commutative algebras first appear in the paper of G.~Baxter~\cite{Baxter}.
For basic results and the main properties of Rota--Baxter algebras see~\cite{GuoMonograph}.

\begin{definition}
Let $A$ be an algebra over a field~$\Bbbk$. A linear operator $R$ on $A$ is called an algebra Rota--Baxter operator (A--RB--operator) of weight~$\lambda\in \Bbbk$ if
\begin{equation}\label{RBAlgebra}
R(a)R(b) = R\left( R(a) b + a R(b) + \lambda a b \right)
\end{equation}
for all $a, b \in A$.
\end{definition}

\begin{example}
In \cite{Ma-1} (see also \cite{Ma})
the algebra RB-operators on $H_4$ as on an associative algebra were found.
It was proven that any such operator has one of the following form:

\begin{itemize}
\item[(a)] $R(1) = 0$,  $R(g) = 0$, $R(x) = -\lambda x$, $R(gx) = -\lambda gx$;

\item[(b)] $R(1) = -\lambda 1$,  $R(g) = -\lambda g$, $R(x) = 0$, $R(gx) = 0$;

\item[(c)] $R(1) = -\lambda 1$,  $R(g) = -\lambda g$, $R(x) = -\lambda x$, $R(gx) = -\lambda gx$;

\item[(d)] $R(1) = 0$, $R(g) = -p_1 1 + p_1 g - \frac{(\lambda + p_1) (\lambda + p_1 + p_2)}{p_3} x + \frac{(\lambda + p_1) (\lambda + p_2)}{p_3} gx$,\\
 $R(x) = -p_3 1 + p_3 g - (2 \lambda + p_1 + p_3) x +  (\lambda + p_2) g x$,\\
   $R(g x) = -p_3 1 + p_3 g - (\lambda + p_1 + p_2) x +  p_2 g x$;

\item[(e)] $R(1) = -\lambda 1$, $R(g) = (\lambda + p_1) 1 + p_1 g - \frac{(\lambda + p_1) (\lambda + p_1 + p_2)}{p_3} x + \frac{(\lambda + p_1) (\lambda + p_2)}{p_3} gx$,\\
 $R(x) = p_3 1 + p_3 g - (2 \lambda + p_1 + p_2) x +  (\lambda + p_2) g x$,\\
   $R(g x) = p_3 1 + p_3 g - (\lambda + p_1 + p_2) x +  p_2 g x$;

\item[(f)] $R(1) = -\lambda 1$, $R(g) = \lambda  1 + p_1  x + \frac{ p_1 p_2}{\lambda + p_2} gx$,\\
$R(x) = - (\lambda +  p_2) x - p_2 g x$, $R(gx) = (\lambda +  p_2) x + p_2 g x$;

\item[(g)] $R(1) = -\lambda 1$, $R(g) = \lambda 1  + \frac{\lambda (\lambda + p_1)}{p_2} x + \frac{\lambda (\lambda + p_1)}{p_2} gx$,\\
 $R(x) = -p_2 1 - p_2 g - (2 \lambda + p_1) x -  (\lambda + p_1) g x$,\\
   $R(g x) = p_2 1 + p_2 g + (\lambda + p_1) x +  p_1 g x$;

\item[(h)]   $R(1) = \frac{1}{2} \lambda  1 -  \frac{1}{2} \lambda g  + p_1 x +  p_2 g x$, $R(g) = \frac{1}{2} \lambda  1 -  \frac{1}{2} \lambda g - p_2 x +  p_1 g x$,\\
$R(x) = -\frac{1}{2} \lambda  x -  \frac{1}{2} \lambda g x$, $R(g x) = -\frac{1}{2} \lambda  x -  \frac{1}{2} \lambda g x$.
\end{itemize}
\end{example}

\medskip

Recall (see \cite{Guo2020}) the definition of group RB-operator.

\begin{definition}
Let $G$ be a group.
\begin{enumerate}
\item A map $B \colon G\to G$ is called a {\it Rota--Baxter operator} of weight~1 if
\begin{equation}\label{RB}
B(g)B(h) = B( g B(g) h B(g)^{-1} )
\end{equation}
for all $g,h\in G$.
\item A map $C \colon G\to G$ is called a {\it Rota--Baxter operator} of weight~$-1$ if
$$
C(g) C(h) = C( C(g) h C(g)^{-1} g )
$$
for all $g,h\in G$.

\end{enumerate}
\end{definition}

There is a~bijection between Rota--Baxter operators of weights~1 and $-1$ on a~group~$G$. We will call group Rota--Baxter
operators of weight~1 simply group Rota--Baxter operators (G--RB--operators).
A group endowed with a Rota--Baxter operator is called a~{\it Rota--Baxter group} (RB--group).

\begin{proposition}[\cite{Guo2020}]\label{prop:Derived}
Let $(G, \cdot ,B)$ be a Rota--Baxter group.

a) The pair $(G, *)$, with the multiplication
\begin{equation}\label{R-product}
g*h = gB(g)hB(g)^{-1},
\end{equation}
where $g,h\in G$, is also a group.

b) The operator $B$ is a Rota--Baxter operator on the group $(G, *)$.

c) The map $B\colon(G, *)\to (G, \cdot)$ is a homomorphism of Rota--Baxter groups.
\end{proposition}

\subsection{Hopf Rota--Baxter operators  on cocommutative Hopf algebras } \label{RRBO}
Let $H$ be a cocommutative Hopf algebra. M.~Goncharov \cite{Goncharov2020} defined an operator $B \colon H \to H$, which is a coalgebra homomorphism and satisfies the identity
\begin{equation} \label{group}
B(a) B(b) = B\left( a_{(1)} \, B(a_{(2)}) \, b \, S(B(a_{(3)})) \right),~~a, b \in H,
\end{equation}
where
$$
\Delta(a) = a_{(1)} \otimes a_{(2)}, ~~\Delta^2(a) = (\id \otimes \Delta ) \Delta(a) = a_{(1)} \otimes a_{(2)} \otimes a_{(3)}.
$$
In particular, if $a$ is a group-like element, then $a = a_{(1)}  = a_{(2)}= a_{(3)}$, and $S(a) = a^{-1}$. In this case we have
$$
B(a) B(b) = B\left( a \, B(a) \, b \, B(a)^{-1} \right).
$$
This is the definition of a Rota--Baxter operator on groups.

We will call a linear operator $B \colon H\to H$ which satisfies the equation (\ref{group}), a \emph{cocommutative Hopf Rota--Baxter operator (CH--RB--operator)} or simply
\emph{Hopf Rota--Baxter operator (H--RB--operator)}.


\subsection{Rota--Baxter operators on arbitrary Hopf algebras}
The next definition is a generalization of the definition of a RB--operator on cocommutative Hopf algebras and it was introduced  in \cite{BN}.

\begin{definition} \label{Def-RRBO}
Let $(H, *, \Delta_H, \varepsilon_H, S_H)$ be a Hopf algebra and the next conditions hold:

1) A map
$$
B \colon H \to H^{op}
$$
is a coalgebra homomorphism.

2) For any $a, b \in H$ holds
\begin{equation}
S\left( B(a_{(2)}) \right) \, b \, B(a_{(3)}) \otimes S\left( B(a_{(1)}) \right) = S\left(B(a_{(2)}) \right) \, b  \, B(a_{(3)}) \otimes S(a_{(1)}) \, a_{(4)} \, S\left( B(a_{(5)})\right).
\end{equation}

3) For any $a, b \in H$ holds
$$
B(a) B(b) = B\left(a_{(1)}  S(B(a_{(2)})) \,  b B(a_{(3)})\right).
$$
We shall call such an operator $B$ a {\it non--cocommutative Hopf Rota--Baxter operator (NCH--RB--operator)} on $H$.
\end{definition}

\bigskip


\section{Connection between group RB--operators and algebra RB--operators} \label{G-A-RBO}

In this section we are studying the following question. Let $A=\Bbbk[G]$ be a  group algebra, then it is an associative algebra and we can find algebra RB--operators on $A$. On the other side, we can find group RB--operators on $G$ and extend them by linearity to $A$. Is there some connection between them?

Let $B$ be a group RB--operator on $G$. Extend it by linearity to the group algebra $\Bbbk[G]$. If the extended operator is a Rota--Baxter operator of weight $\lambda$ on $\Bbbk[G]$, i.e.
$$
B(g)B(h)=B(B(g)h+gB(h)+\lambda gh),\quad g,h\in G,
$$
we say that $B$ is a group-algebra Rota--Baxter operator (GA--RB--operators) of weight $\lambda$.

\begin{proposition}
Let $B$ be a group RB--operator on $G$. Then $B$ is  an algebra RB--operator of weight $\lambda$ on $\Bbbk[G]$ if and only if $\lambda=-1$ and for any $g,h\in G$ we have either
\begin{equation}\label{eq:garbo}
\left\{
\begin{array}{l}
     B(g)B(h)=B(B(g)h),  \\
     B(gB(h))=B(gh),
\end{array}
\right.
\qquad\mbox{or}\qquad
\left\{
\begin{array}{l}
     B(g)B(h)=B(gB(h)),  \\
     B(B(g)h)=B(gh).
\end{array}
\right.
\end{equation}
\end{proposition}

\begin{proof}
Suppose that $B \colon G \to G$ is a G--RB--operator, then $B(g) \in G$ and $B$ satisfies the identity
$$
B(g)B(h)=B(B(g)h)+B(gB(h))+\lambda B(gh),\quad g,h\in G,
$$
if the  left hand side is equal to one summand from the right hand side and sum of the other two terms is zero. But $B(B(g)h)+B(gB(h)) \not = 0$. Hence,
$$
B(B(g)h) = - \lambda B(gh) ~\mbox{or}~B(gB(h)) = - \lambda B(gh).
$$
It is possible if $\lambda=-1$ and summands satisfies one condition from proposition.
\end{proof}

\begin{proposition}
    Let $B$ be a GA--RB--operator on $G$. Then $B^2=B$.
\end{proposition}
\begin{proof}
    Let $h=1$. By~\eqref{eq:garbo} either $B(g)=B(g)B(1)=B(B(g)1)=B^2(g)$ or $B^2(g)=B(B(g)1)=B(g1)=B(g)$.
\end{proof}

Let $B$ be a group RB--operator. Denote $H=\im B$ and $K=\ker B$. Then $H$ and $K$ are subgroups of $G$ (see \cite{Guo2020}).

\begin{theorem} \label{GA-oper}
Let $B$ a GA--RB--operator on $G$. Then
\begin{enumerate}
    \item $G=KH$, where  $H=\im B$ and $K=\ker B$, is an exact factorization;
    \item $B(kh)=h$ for all $k \in K$, $h\in H$;
    \item $H$ is commutative and $h^2Kh^{-2}\subset K$ for all $h\in H$.
\end{enumerate}

Conversely, for any exact group factorization $G=KH$ such that $H$ is commutative and $h^2Kh^{-2}\subset K$ for all $h\in H$, the operator $B(kh)=h$ is a GA--RB--operator on $G$.
\end{theorem}

\begin{proof}
  1.  Let $B$ be a GA--RB--operator. For any $g\in G$ take $h=B(g)^{-1}$. Then $h\in H$ and $B(h)=h$. Hence
$$
B(gB(g)^{-1})=B(gB(g)B(g)^{-1}B(g)^{-1})=B(g)B(B(g)^{-1})=B(g)B(g)^{-1}=1.
$$
Then $g=(gB(g)^{-1})B(g)\in KH$.

Let $g\in K\cap H$. Then $g=B(g)=1$. Thus, $G=KH$ is an exact factorization.

For any $g=kh\in G$, $B(g)=B(kh)=B(kB(k)hB(k)^{-1})=B(k)B(h)=B(h)=h$.

As a consequence, we have $B(gh)=B(g)h$ for any $g\in G$, $h\in H$.

Let $h_1,h_2\in H$. Then $h_1h_2=B(h_1)B(h_2)=B(h_1h_1h_2h_1^{-1})=h_1^2h_2h_1^{-1}$. Hence, $h_1h_2=h_2h_1$ and $H$ is commutative.

Let $k\in K$ and $h\in H$. Then $1=B(h)B(k)h^{-1}=B(h^2kh^{-1})h^{-1}=B(h^2kh^{-2})$. Thus, $h^2kh^{-2}\in K$.

2. Now, let $G=KH$ be an exact factorization such that $H$ is commutative and $h^2Kh^{-2}\subset K$ for all $h\in H$, and $B(kh)=h$, $k\in K, h\in H$. Then for any $g\in G, k\in K, h\in H$  $B(kgh)=B(g)h$.

Let $g_1=k_1h_1, g_2=k_2h_2$. Check the conditions~\eqref{eq:garbo}:
$$B(g_1)B(g_2)=h_1h_2=B(k_1h_1h_2)=B(g_1B(g_2)),$$ $$B(B(g_1)g_2)=B(h_1k_2h_2)=B(k_1h_1k_2h_2)=B(g_1g_2).$$

Let us check the group RB--operator condition
\begin{multline*}
B(g_1)B(g_2)=h_1h_2=B(h_1^2k_2h_1^{-2})h_1h_2=B(h_1^2k_2h_1^{-1}h_2)=B(k_1h_1^2k_2h_2h_1^{-1})=\\ B(g_1B(g_1)g_2B(g_1)^{-1}).
\end{multline*}
Here we used the condition $h_1^2k_2h_1^{-2}\in K$ and $h_1^{-1}h_2=h_2h_1^{-1}$.
Thus, $B$ is a GA--RB--operator.
\end{proof}

\bigskip


\section{Hopf and algebra Rota--Baxter operators on $H_4$} \label{RB_H4}

In the present section we find cocommutative Hopf RB--operators (CH--RB--operator) $B \colon H_4 \to H_4$, that means that $B$ satisfy the identity
\begin{equation} \label{FocCocom}
B(a) B(b) = B \left( a_{(1)} \, B(a_{(2)}) \, b \, S(B(a_{(3)})) \right),~~a, b \in H,
\end{equation}
where
$$
\Delta^2(a) = (\Delta \otimes \id ) \Delta (a) =  a_{(1)} \otimes  a_{(2)} \otimes a_{(3)},
$$
which are extensions of group RB--operators. Let  $G$ be  a cyclic group of order 2 group, which lies in $H_4$ and is generated by $g$. On this group there  exist two different RB--operators:
$B_0$ and $B_{-1}$ such that
$$
B_0(1) = 1,~~B_0(g) = 1;~~~B_{-1}(1) = 1,~~B_{-1}(g) = g.
$$
By linearity these operators can be extended to group algebra $\Bbbk[G]$. We have to extend these operators on $H_4$. In order to do this, we must define these operators on $x$ and $gx$.

At first, let us extend $B_0$. We put
$$
B_0(x) = a_1 + a_2 g + a_3 x + a_4 gx,~~~B_0(gx) = b_1 + b_2 g + b_3 x + b_4 gx,~~~a_i, b_i \in \Bbbk.
$$

\begin{theorem} \label{t4.1}
\begin{enumerate}
\item If $char(\Bbbk) = 2$, then the map
$$
B_0(1) =B_0(g) = 1,~~ B_0(x) = B_0(gx) = a_1 + a_1g + a_3 x + a_4 gx,~~a_i \in  \Bbbk,
$$
is a CH--RB--operator  on $H_4$.

\item If $char(\Bbbk) \not= 2$, then

a) $B_0(1) =B_0(g) = 1,~~ B_0(x) = B_0(gx) = a_1 + a_1g + a_3 x - 2 a_3 gx,$

and

b) $B_0(1) =B_0(g) = 1,~~ B_0(x) = B_0(gx) = a_1 - a_1g + a_3 x - 2 a_3 gx,~~a_i \in  \Bbbk$

are  CH--RB--operators  on $H_4$.
\end{enumerate}
\end{theorem}

\begin{proof}

It is necessary to find conditions under which equality (\ref{FocCocom}) holds. It is enough to check it for the basic elements.

1) If we put $a=g$ and $b=x$ in (\ref{FocCocom}),  we get
$$
B_0 (g) B_0(x) = B_0\left( g \, B_0(g) \, x \, S(B_0(g))\right).
$$
By definition $B_0$ and $S$ on $g$, we get
$$
B_0(x) = B_0(g x).
$$

2) Further, if we put $a=x$ and $b=g$ in (\ref{FocCocom}),  we get
$$
B_0 (x) B_0(g) = B_0\left( x_{(1)} \, B_0( x_{(2)}) \,  g \, S(B_0( x_{(3)}))\right). 
$$
Using the equality
$$
\Delta^2(x) = (\Delta \otimes \id ) \Delta (x) =  (\Delta \otimes \id ) ( x \otimes 1 + g \otimes x) =  x \otimes 1 \otimes 1 + g \otimes x \otimes 1 + g \otimes g \otimes x,
$$
we have
$$
B_0(x) = B_0(x g + g B_0(x) g + S(B_0(x))).
$$
By 1), $B_0(x g) = -B_0(x)$. Hence,
$$
B_0(x) = -B_0(x) + B_0(g B_0(x) g) + B_0(S(B_0(x))) ~\Leftrightarrow ~ 2 B_0(x) =  B_0(g B_0(x) g) + B_0(S(B_0(x))).
$$
Since
$$
B_0(g \, B_0(x) \,  g) = B_0(a_1 + a_2 g - a_3 x - a_4 gx),~~ B_0(S(B_0( x))) = B_0(a_1 + a_2 g - a_3 g x + a_4 x),
$$
we have the equality
$$
2 B_0 (a_1 + a_2 g + a_3 x + a_4 gx)  = B_0(a_1 + a_2 g - a_3 x - a_4 gx) +  B_0(a_1 + a_2 g - a_3 g x + a_4 x),
$$
which is  equivalent to
$$
 B_0((3 a_3 - a_4) x + (3 a_4 + a_3) gx) = 0 ~\Leftrightarrow ~ (4 a_3 + 2a_4) B_0(x) = 0.
$$

If characteristic $\Bbbk$ is 2, then $B_0(x) = B_0(g x)$ can be arbitrary element in $H_4$.

If $char(\Bbbk) \not= 2$, then $B_0(x)= B_0(g x) = 0$ or
$$
 B_0(x) = B_0(gx) = a_1 + a_2g + a_3 x - 2 a_3 gx.
$$

3) If we put $a=b=x$ in (\ref{FocCocom}),  we get
$$
B_0 (x) B_0(x) = B_0\left( x_{(1)} \, B_0( x_{(2)}) \,  x \, S(B_0( x_{(3)}))\right).
$$
It is equivalent to
$$
B_0 (x)^2 = B_0\left(g B_0(x)  x + g x S(B_0(x))\right).
$$
Since
$$
g B_0(x) x = a_1 g x + a_2 x ,~~~g x S(B_0(x)) = a_1 gx - a_2 x,
$$
we have
$a_1^2  - a_2^2 = 0$. Hence,  if $char(\Bbbk) = 2$, then $a_1 = a_2$. If $char(\Bbbk) \not= 2$, then $a_2 = a_1$ or $a_2 = -a_1$.

4) If we put $a=b=gx$ in (\ref{FocCocom}),  we get
$$
B_0 (g x) B_0(g x) = B_0\left( (gx)_{(1)} \, B_0( (gx)_{(2)}) \,  g x \, S(B_0( (gx)_{(3)}))\right).
$$
Using the equality
$$
\Delta^2(gx) = (\Delta \otimes \id ) \Delta (gx) =  (\Delta \otimes \id ) ( g x \otimes g + 1 \otimes g x) =  g x \otimes g \otimes g + 1 \otimes g x \otimes g + 1 \otimes 1 \otimes g x,
$$
we have
$$
B_0 (g x) B_0(g x) = B_0\left( B_0(g x) g x + g x S(B_0(gx)\right).
$$
Since $B_0 (g x) = B_0 (x)$, we get
$$
B_0 (x)^2 = B_0\left( B_0(x) g x + g x S(B_0(x)\right).
$$
This is equality from 3).

By inserting in (\ref{FocCocom}) other basic elements of $H_4$ we can check that we do not get other restrictions on $B_0$.

\end{proof}

Now let us consider extensions of $B_{-1}$ on $H_4$.
Put
$$
B_{-1}(x) = a_1 + a_2 g + a_3 x + a_4 gx,~~~B_{-1}(gx) = b_1 + b_2 g + b_3 x + b_4 gx,~~~a_i, b_i \in \Bbbk.
$$
The next theorem holds

\begin{theorem} \label{t4.2}
\begin{enumerate}
\item If $char(\Bbbk) = 2$, then there exists a CH--RB--operator on $H_4$:

$B_{-1}(1) =1,~~B_{-1}(g) = g,~~ B_{-1}(x) = a_1 + a_1 g + a_3 x + a_4 gx,~~ B_{-1}(gx) = a_1 + a_1 g + a_4 x + a_3 gx,~~a_i \in  \Bbbk.$\\

\item If $char(\Bbbk) \not= 2$, then there exist three  types of CH--RB--operators on $H_4$:

a) $B_{-1}(1) =1,~~B_{-1}(g) = g,~~  B_{-1}(x) = B_{-1}(gx) =  a - a g,~~$  $a \in \Bbbk$;

b) $B_{-1}(1) =1,~~B_{-1}(g) = g,~~  B_{-1}(x) = B_{-1}(gx) = - x + g x$;

c)  $B_{-1}(1) =1,~~B_{-1}(g) = g,~~  B_{-1}(x) = - B_{-1}(gx) = a + a g + x + g x,$  $a \in \Bbbk$.
\end{enumerate}
\end{theorem}

\begin{proof}

It is necessary to find conditions under which equality (\ref{FocCocom}) holds. It is enough to check it for the basic elements.

1) If we put $a=g$ and $b=x$ in (\ref{FocCocom}),  we get
$$
B_{-1} (g) B_{-1}(x) = B_{-1}\left( g \, B_{-1}(g) \, x \, S(B_{-1}(g))\right).
$$
From this identity, by definition $B_{-1}$ and $S$ on $g$, we get
$$
g B_{-1}(x) = - B_{-1}(g x).
$$
Hence,
$$
B_{-1}(g x) = - a_2 - a_1 g - a_4 x - a_3 gx.
$$

2)  If we put $a=x$ and $b=g$ in (\ref{FocCocom}),  we get
$$
B_{-1} (x) B_{-1}(g) = B_{-1}\left( x_{(1)} \, B_{-1}( x_{(2)}) \,  g \, S(B_{-1}( x_{(3)}))\right)
$$
or
$$
B_{-1} (x) g = B_{-1}\left( x g + g B_{-1}(x) g + g S(B_{-1}(x)) \right).
$$
Using the formulas
$$
B_{-1} (x g) =  - B_{-1} (g x) =  g B_{-1}(x),
$$
we can find terms from the right hand side:
$$
 B_{-1} (g B_{-1}(x) g) = B_{-1}(a_1 +a_2 g - a_3 x - a_4 gx) = a_1 +a_2 g + (- a_3 + a_4 g) B_{-1}(x),
$$
$$
  B_{-1}(g S(B_{-1}(x))) = B_{-1}(a_1g  +a_2 - a_3 x + a_4 gx) = a_1 g +a_2  - a_3  B_{-1}(x)-  a_4 g  B_{-1}(x) =
$$
$$
=  a_2 +a_1 g  - (a_3 + a_4 g)  B_{-1}(x),
$$
and we get
$$
B_{-1} (x) g = (g - 2 a_3) B_{-1}(x) + a_1 +a_2 + (a_1 +a_2) g.
$$
Inserting the equality for $B(x)$, we have
$$
0 = (a_1 +a_2 - 2 a_1 a_3)  +  (a_1 +a_2 - 2 a_2 a_3) g + 2(a_4 - a_3^2) x + 2(a_3 - a_3 a_4) g x.
$$
From this equality we get to the system
$$
\begin{cases}
a_1 + a_2 - 2 a_1 a_3 = 0,  \\
a_1 + a_2 - 2 a_2 a_3 = 0,  \\
2 (a_4 - a_3^2) = 0,  \\
2 (a_3 - a_3 a_4) = 0.  \\
\end{cases}
$$

If $char(\Bbbk) = 2$, then this system is equivalent to the equality $a_2 = a_1$.This is the map from (1).

Suppose that $char(\Bbbk) \not= 2$, then from the last equation follows that $a_3 = 0$, or $a_3 \not= 0$ and $a_4 = 1$.

If $a_3 = 0$, then the system is equivalent to the system
$$
\begin{cases}
a_1 + a_2  = 0,  \\
 a_4 = 0.  \\
\end{cases}
$$
Its solutions $(a_1, a_2, a_3, a_4) = (a_1, -a_1, 0, 0)$, $a_1 \in \Bbbk$. This is the map from (2) case a).

If $a_3 \not = 0$ and $a_4 = 1$,  then the system is equivalent to the system
$$
\begin{cases}
a_1 + a_2 - 2 a_1 a_3 = 0,  \\
a_1 + a_2 - 2 a_2 a_3 = 0,  \\
 a_3^2 = 1.  \\
\end{cases}
$$
Its solutions  $(a_1, a_2, a_3, a_4) = (0, 0, -1, 1)$ and $(a_1, a_2, a_3, a_4) = (a_1, a_1, 1, 1)$, $a_1 \in \Bbbk$. These maps are maps from (2) cases b) and c).

We have three maps which are defined by guadruple
$$
(a_1, a_2, a_3, a_4) = (a, -a, 0, 0),~~(a_1, a_2, a_3, a_4) = (0, 0, -1, 1),~~(a_1, a_2, a_3, a_4) = (a_1, a_1, 1, 1)
$$
and they are CH--RB--operators for some basis elements. It is need to check that these maps are CH--RB--operators for other basis elements.

Let us consider the first guadruple  $(a, -a, 0, 0)$. It defines the  map
$$
B_{-1}(x) = a + a g,~~B_{-1}(g x) = -a - a g.
$$
Let us check that it  is a CH--Rota--Baxter operator for $a = g x$, $b = g$. We must check the equality
$$
B_{-1} (gx) B_{-1}(g) = B_{-1}\left( (gx)_{(1)} \, B_{-1}( (gx)_{(2)}) \,  g \, S(B_{-1}( (gx)_{(3)}))\right).
$$
Since,
$$
\Delta^2(g x) = gx \otimes g \otimes g +  1\otimes gx \otimes g + 1 \otimes 1 \otimes gx,
$$
we get
$$
B_{-1} (gx) B_{-1}(g) = B_{-1}\left( gx \, B_{-1}(g) \,  g \, S(B_{-1}(g)) +   B_{-1}(gx) \,  g \, S(B_{-1}(g)) + g \, S(B_{-1}(gx))  \right).
$$
This equality is equivalent to
\begin{equation} \label{eqRB}
B_{-1} (gx) g = B_{-1}\left( -x +   B_{-1}(gx) + g \, S(B_{-1}(gx))  \right).
\end{equation}
Since $B_{-1} (gx)  = - g B_{-1}(x)$, we have
$$
B_{-1} (gx) + g \, S(B_{-1}(gx)) = - 2 a_2 - 2 a_1 g - (a_3 + a_4) x + (- a_3 + a_4) g x.
$$
Hence,  (\ref{eqRB}) has the form
$$
- g B_{-1} (x) g =  -B_{-1}(x) -  2 a_2 - 2 a_1 g - (a_3 + a_4) B_{-1}(x) + ( a_3 - a_4) g B_{-1}(x).
$$
Take into account that $a_2 = -a_1$ and $a_3 = a_4 = 0$ we have
$$
 g B_{-1} (x) g =  B_{-1}(x).
$$
 It easy to check that this equality is true.

By analogy we can check the map $B_{-1}$ which is defined by quadruple  $(a, -a, 0, 0)$ is a CH--RB--operator for other basis elements.

Further we have to consider two other quadruples: $(0,0, -1, 1)$ and $(a, a, 1, 1)$. Analysing Rota--Baxter relation for basis elements, we see that we have CH--Rota--Baxter operators.

\end{proof}

It is known \cite{Goncharov2020} that a group RB--operator  on a cocommutative Hopf algebra is a cohomomorphism,  that means
$$
(B \otimes B) \Delta = \Delta B.
$$
Let us check this equality for operator $B = B_{-1}$ in Theorem \ref{t4.2} (2) a). Consider the equality
$$
(B \otimes B) \Delta (x)= \Delta B (x).
$$
Its left hand side
$$
(B \otimes B) \Delta (x)= (B \otimes B) (x \otimes 1 + g \otimes x) = (a - ag) \otimes 1 + g \otimes (a - ag) =
$$
$$
= a (1 \otimes 1) - a  (g \otimes g).
$$
Its right hand side
$$
\Delta B (x) = \Delta (a - ag) = a (1 \otimes 1 ) - a  (g \otimes g).
$$
Hence, the equality
$$
(B \otimes B) \Delta (x)= \Delta B (x).
$$
holds.

Also, we can check the equality
$$
(B \otimes B) \Delta (x)= \Delta B (x).
$$

\begin{proposition} \label{notHom}
There exist CH--RB--operator on $H_4$ which is not a coalgebra homomorphism.
\end{proposition}

\begin{proof}
Consider a CH-RB-operator $B_0$ in Theorem 4.1 such that $a_1\ne 0$. Then $\Delta B(x)$ contains the monomial $g\otimes g$ with the coefficient $a_1$. On the other hand, the expression $(B\otimes B)\Delta(x)$ does not contain this monomial. Hence, $\Delta B\ne (B\otimes B)\Delta$.
\end{proof}

As we know the Hopf algebra $H_4$ is not cocommutative.

\subsection{Comparing of CH--RB--operators on  $H_4$ and LA--RB--operators on $H_4^{(-)}$}
In the previous subsection we constructed CH--RB--operators on  $H_4$.
Denote by $L = H_4^{(-)}$ the Lie algebra on the associative algebra $H_4$ in which multiplication  is defined by the rule $[a, b] = ab - ba$ for $a, b \in H_4$.
All Lie algebra RB--operators on $L$ over a field of $char(\Bbbk) \not= 2$ were found in \cite{BB1}.

Further we assume that $char(\Bbbk) \not= 2$.


Let us check what of CH--RB--operators on $H_4$ are Lie algebra RB--operators on  $L$, i.~e.
 satisfies the equality
\begin{equation} \label{LRB}
[B(a), B(b)] = B ([a, B(b)] + [B(a), b] + \lambda [a, b]),~~a, b \in L.
\end{equation}
It is easy to see that if we permute $a$ and $b$ in this equality, we get an equivalent equality.

{\it Case 1:}
Let us take the  RB--operator on $H_4$ from Theorem \ref{t4.2} (2) a):
\begin{equation} \label{Lie-1}
B_{-1}(1) = 1,~~B_{-1}(g) = g,~~B_{-1}(x) = B_{-1}(gx) =  a - a g,~~a \in \Bbbk,
\end {equation}
and check: is it  a Lie algebra RB--operator on $L$? It is easy to see that if $a = b$, or $a=1$, or $b=1$ the equality (\ref{LRB}) hold for any $\lambda \in \Bbbk$. Further we will consider other possibilities.

1) If $a = x$, $b = g$, we get
$$
[B_{-1}(x), B_{-1}(g)] = B_{-1} ([x, B_{-1}(g)] + [B_{-1}(x), g] + \lambda [x, g]).
$$
Since,
$$
[B_{-1}(x), B_{-1}(g)] = 0,~~[x, B_{-1}(g)] = [x, g] = -2 g x,~~ [B_{-1}(x), g] = 0,
$$
we have
$$
0 = -2 B_{-1}\left(gx - \lambda gx \right) \Leftrightarrow 0 =  -2(1 + \lambda)  B_{-1}(gx) \Leftrightarrow 0 =  -2 (1 + \lambda)  a (1 - g).
$$
 It is possible if and only if $ \lambda  = -1$ and $a$ is arbitrary, or $a = 0$ and $\lambda$ is arbitrary.

2) If $a = g x$, $b = g$, we get
$$
[B_{-1}(g x), B_{-1}(g)] = B_{-1} ([g x, B_{-1}(g)] + [B_{-1}(g x), g] + \lambda [g x, g]).
$$
Since,
$$
[B_{-1}(gx), B_{-1}(g)] = 0,~~[gx, B_{-1}(g)] = [gx, g] = -2  x,~~ [B_{-1}(gx), g] = 0,
$$
we have
$$
0 = B_{-1}(-2x - 2 \lambda x) \Leftrightarrow  0 = -2 (1 + \lambda)  a (1 - g).
$$
As in 1) it is possible if and only if $ \lambda  = -1$ and $a$ is arbitrary, or $a = 0$ and $\lambda$ is arbitrary.

Check other relations, we get the same conditions.

{\it Case 2:}
Let us take the  RB-operator on $H_4$  from   Theorem \ref{t4.2} (2) b):
\begin{equation} \label{Lie-3}
B_{-1}(1) = 1,~~B_{-1}(g) = g,~~B_{-1}(x) =  B_{-1}(gx) = - x + g x,
\end {equation}
and find conditions under which  it is an algebra RB--operator on $L$.

If $a =  x$, $b = g$, we get
$$
[B_{-1}(x), B_{-1}(g)] = B_{-1} ([x, B_{-1}(g)] + [B_{-1}(x), g] + \lambda [x, g]).
$$
Let us find the products:
$$
[B_{-1}(x), B_{-1}(g)] = [B_{-1}(x), g] = [- x + g x, g] = 2( -x + gx),
$$
$$
[x, B_{-1}(g)] = [x, g] = -2gx.
$$
Then
$$
2(x - g x)  = 2 B_{-1}(x +  \lambda g x) \Leftrightarrow (2 + \lambda) x - (2 + \lambda) g x = 0.
$$
From this relation follows that  $ \lambda  =-2$.

Check other relations, we prove that the operator, constructed in (\ref{Lie-3}) is a Lie algebra RB--operator on $L$ of weight $-2$.

{\it Case 3:}
Let us take the  RB--operator on $H_4$ from   Theorem \ref{t4.2} (2) c):
\begin{equation} \label{Lie-2}
B_{-1}(1) = 1,~~B_{-1}(g) = g,~~B_{-1}(x) = - B_{-1}(gx) = a + a g + x + g x,~~a \in \Bbbk,
\end {equation}
and check that it is a Lie algebra RB--operator on $L$.

If $a =  x$, $b = g$, we get
$$
[B_{-1}(x), B_{-1}(g)] = B_{-1} \left([x, B_{-1}(g)] + [B_{-1}(x), g] + \lambda [x, g]\right).
$$
Let us find the products:
$$
[B_{-1}(x), B_{-1}(g)] = [a + a g + x + g x, g] = -2(x + gx),
$$
$$
[x, B_{-1}(g)] = [x, g] = -2gx,~~~[B_{-1}(x), g] =  -2(x + gx).
$$
Then
$$
-2(x + gx) = -2 B_{-1}(g x + \lambda gx).
$$
Take into account the values of $B_{-1}$, we get
$$
(1 + \lambda) a + (1 + \lambda) a g + (2 + \lambda) x + (2 + \lambda) g x = 0.
$$
From this relation follows that  $a= 0$ and $\lambda  =-2$.

Check other relations, we prove that the operator, constructed in (\ref{Lie-2}) is a Lie algebra RB--operator on $L$ if $a= 0$ and $\lambda  =-2$.

Hence, we get

\begin{proposition} Let $char(\Bbbk) \not= 2$.
\begin{enumerate}
\item For any  $ \lambda\in\Bbbk$ the map
$$
B_{-1}(1) = 1,~~B_{-1}(g) = g,~~B_{-1}(x) =  B_{-1}(gx) = 0,
$$
defines a Lie algebra RB--operator of weight $\lambda$ on $L$.
\item The map
$$
B_{-1}(1) = 1,~~B_{-1}(g) = g,~~B_{-1}(x) = B_{-1}(gx) =  a - a g,~~a \in \Bbbk,
$$
defines a Lie algebra RB--operator of weight $ \lambda  =-1$ on $L$.
\item The map
$$
B_{-1}(1) = 1,~~B_{-1}(g) = g,~~B_{-1}(x) =  B_{-1}(gx) = - x + g x,
$$
defines a Lie algebra RB--operator of weight $ \lambda  =-2$ on $L$.
\item The map
$$
B_{-1}(1) = 1,~~B_{-1}(g) = g,~~B_{-1}(x) = - B_{-1}(gx) =  x + g x,
$$
defines a Lie algebra RB--operator of weight $ \lambda  =-2$ on $L$.
\end{enumerate}
\end{proposition}

Let us consider CH--RB--operators  $B_0$, which were constructed in Theorem \ref{t4.1} and find conditions under which these operators are algebra RB--operators on $L$, that means that for some $\lambda \in \Bbbk$ and any $a, b \in L$ holds
\begin{equation} \label{LRB-1}
[B_0(a), B_0(b)] = B_0 ([a, B_0(b)] + [B_0(a), b] + \lambda [a, b]),~~a, b \in L.
\end{equation}
It is enough to consider the cases: $a, b \in \{ 1, g, x, g x \}$. It is easy to see that if $a, b \in \{ 1, g \}$ the need equality holds for any $\lambda$. As we see before, if the equality holds to a pair $(a, b)$, then it holds for the pair $(b, a)$. Hence, it is need to consider only one possibility: $a = gx$, $b = x$. We will analyze two cases from Theorem \ref{t4.1} (2).

{\it Case 1}: $B_0(1) =B_0(g) = 1,~~ B_0(x) = B_0(gx) = a_1 + a_1g + a_3 x - 2 a_3 gx.$ We have to find conditions under which holds
\begin{equation} \label{LRB-1}
[B_0(g x), B_0(x)] = B_0 \left( [gx, B_0(x)] + [B_0(gx), x] + \lambda [gx, x] \right).
\end{equation}
At first, let us find the products
$$
[B_0(g x), B_0(x)] = 0,~~[gx, B_0(x)] = [gx, a_1 + a_1g + a_3 x - 2 a_3 gx] = -2a_1 x,
$$
$$
[B_0(gx), x] = 2 a_1 g x,~~[gx, x] = 0.
$$
If $a=g$, $b=g$, then $(2\lambda - 2a_3)B_0(x)=0$. Hence, $B_0(x)=0$ or $\lambda=a_3$.

{\it Case 2}: $B_0(1) =B_0(g) = 1,~~ B_0(x) = B_0(gx) = a_1 - a_1g + a_3 x - 2 a_3 gx.$ We have to find conditions under which holds
\begin{equation} \label{LRB-2}
[B_0(g x), B_0(x)] = B_0 \left( [gx, B_0(x)] + [B_0(gx), x] + \lambda [gx, x] \right).
\end{equation}
At first, let us find the products
$$
[B_0(g x), B_0(x)] = 0,~~[gx, B_0(x)] = [gx, a_1 - a_1g + a_3 x - 2 a_3 gx] = 2 a_1 x,
$$
$$
[B_0(gx), x] = -2 a_1 g x,~~[gx, x] = 0.
$$
Then (\ref{LRB-2}) has the form
$$
0 =  2 a_1 B_0( x -  g x).
$$
Since, $B_0( x -  g x) = 0$, the equality (\ref{LRB-2}) holds for any  $\lambda, a_1, a_3 \in \Bbbk$.

Hence, we get

\begin{proposition} Let $char(\Bbbk) \not= 2$.
\begin{enumerate}
\item For any  $ \lambda, a \in \Bbbk$ the map
$$
B_0(1) = B_0(g) = 1,~~B_0(x) =  B_0(gx) = a x - 2 a g x,
$$
defines a Lie algebra RB--operator of weight $\lambda$ on $L$.
\item  For any  $\lambda, a_1, a_3 \in \Bbbk$ the map
$$
B_0(1) =B_0(g) = 1,~~ B_0(x) = B_0(gx) = a_1 - a_1g + a_3 x - 2 a_3 gx
$$
defines a Lie algebra RB--operator of weight $\lambda$ on $L$.

\end{enumerate}
\end{proposition}

In this subsection we assumed that $char(\Bbbk) \not= 2$. What can we say in  the case $char(\Bbbk) = 2$?

\bigskip


\section{Some questions for future research}

\begin{conjecture}
Any Lie algebra RB--operator on  $H^{(-)}_4$ is an associative algebra RB--operator on $H_4$ or a cocommutative Hopf RB--operator on $H_4$.
\end{conjecture}

\begin{question}
Let $H$ be a cocommutative Hopf algebra other a field of characteristic 0 and $B \colon H \to H$ such that

1)  $B(a) B(b) = B \left( a_{(1)} \, B(a_{(2)}) \, b \, S(B(a_{(3)})) \right),~~a, b \in H$,

2) $B$ sends a group like--element to a group--like element,

3) $B$ sends a primitive element to a primitive element.

Is it true that $B$ is a coalgebra homomorphism?
\end{question}

As we proved in Proposition \ref{notHom} for non-cocommutative Hopf algebras the answer to this question is no.

\begin{question}
It is known that any cocommutative Hopf algebra $H$ over a field of characteristic~0 is a semi-direct product $H = \Bbbk[G] \otimes U(L)$ of group algebra $\Bbbk[G]$ and the universal enveloping algebra $U(L)$ of Lie algebra $L$ of primitive elements (Milnor--Moore theorem). The product of two elements in $H$ is defined
$$
(g_1 \otimes l_1) (g_2 \otimes l_2) = (g_1 g_2 \otimes l_1^{g_2} l_2).
$$
Suppose that
$$
B_1 \colon \Bbbk[G]  \to \Bbbk[G],~~B_2 \colon U(L)  \to U(L)
$$
are group and algebra RB--operator, respectively. Under which conditions  $B_1 \otimes B_2$ is a cocommutative Hopf RB--operator on $H$?
\end{question}

\begin{question}
In Theorems \ref{t4.1} and \ref{t4.2} we have found cocommutative Hopf RB--operators on  $H_4$, which are extensions of  group RB--operators on the cyclic group $\{ 1, g \}$. Are there cocommutative Hopf RB--operators on $H_4$ which are not extensions of group RB--operators on this cyclic group?
\end{question}

\begin{question}
Find non--cocommutative Hopf RB--operators (in the sense \cite{BN}) on $H_4$.
\end{question}

\begin{question}
Find non--cocommutative Hopf RB--operator (in the sense \cite{BN})  on $U_q(sl_2)$.
\end{question}

\begin{question}
M. Goncharov \cite{Goncharov2020} considered RB--operators on Lie algebra $sl_2$ We know RB--operator on $U(sl_2)$ and extended them to  RB--operator on $U(sl_2)$ as on Hopf algebra. What will we get if take the Lie algebra of $U(sl_2)$ and find Lie algebra RB--operators on this algebra?
\end{question}

\begin{question} (P. Kolesnikov.)
The algebra $U_q(sl_2)$ is a deformation of the algebra $U(sl_2)$. Is it possible, using this deformation, construct RB--operators on $U_q(sl_2)$ if we know RB--operators on $U(sl_2)$?
\end{question}

\bigskip


\section*{Acknowledgments}
Authors are grateful to  V.~Gubarev, M.~Goncharov, A. Pozhidaev,  and P. Kolesnikov for the fruitful discussions and useful suggestions.
Authors are also grateful to participants of the seminar ``\'{E}variste Galois''
at Novosibirsk State University for attention to our work.

Valery G. Bardakov is supported by the Russian Science Foundation (RSF
 24-21-00102) for work in section 3.

Viktor N. Zhelyabin is supported by the Program of Fundamental Research RAS,
	(FWNF-2022-0002) for work in section  4.


\medskip

\noindent Valeriy G. Bardakov \\
Sobolev Institute of Mathematics,  Acad. Koptyug ave. 4, 630090 Novosibirsk, Russia; \\
Novosibirsk State Agrarian University,
Dobrolyubova str., 160, 630039 Novosibirsk; \\
Regional Scientific and Educational Mathematical Center of Tomsk State University; \\
Lenin ave. 36, 634009 Tomsk, Russia; \\
email: bardakov@math.nsc.ru

\medskip
\noindent  Viktor N. Zhelaybin \\
Sobolev Institute of Mathematics,  Acad. Koptyug ave. 4, 630090 Novosibirsk, Russia; \\
email: vicnic@mail.math.nsc.ru

\medskip
\noindent Igor Nikonov \\
Lomonosov Moskow  State University; \\
email: nikonov@mech.math.msu.su

\end{document}